\documentclass{amsart}
\usepackage[headings]{fullpage}
\usepackage{amssymb,epic,eepic,epsfig,amsbsy,amsmath,amscd}

\theoremstyle{plain}

\newtheorem{theorem}{Theorem}
\newtheorem{proposition}{Proposition}[section]
\newtheorem{lemma}[proposition]{Lemma}
\newtheorem{corollary}[proposition]{Corollary}

\theoremstyle{definition}

\theoremstyle{remark}

\newtheorem{remark}[proposition]{Remark}

\def\printname#1{
    \if\draft y
        \smash{\makebox[0pt]{\hspace{-0.5in}
            \raisebox{8pt}{\tt\tiny #1}}}
    \fi
}

\newlength{\standardunitlength}
\catcode`\@=11
\long\def\@makecaption#1#2{%
     \vskip 10pt

\setbox\@tempboxa\hbox{
       \small\sf{\bfcaptionfont #1. }\ignorespaces #2}%
     \ifdim \wd\@tempboxa >\captionwidth {%
         \rightskip=\@captionmargin\leftskip=\@captionmargin
         \unhbox\@tempboxa\par}%
       \else
         \hbox to\hsize{\hfil\box\@tempboxa\hfil}%
     \fi}
\font\bfcaptionfont=cmssbx10 scaled \magstephalf
\newdimen\@captionmargin\@captionmargin=2\parindent
\newdimen\captionwidth\captionwidth=\hsize
\catcode`\@=12

\newcommand{\qbinom}[2]{\text{$\left[\begin{array}{c}#1\\ #2\end{array}
\right]$}}

\def\BN{\mathbb N}
\def\BZ{\mathbb Z}

\def\BQ{\mathbb Q}

\def\BC{\mathbb C}

\def\cL{\mathcal L}

\newcommand{\sn}{\operatorname{sn}}

\def\Z{\BZ}

\def\msl{\mathfrak{sl}}

\def\ev{\mathrm{ev}}

\newcommand{\Gd} {\Lambda_d}
\newcommand{\hGd} {\hat \Lambda_d}
\newcommand{\Ud} {U_d}

\newcommand{\Sd} {\BN_d}

\begin{document}


\title[Quantum invariants of 3-manifolds]{Strong Integrality of Quantum Invariants of 3-manifolds}

\author{Thang T. Q. L\^e}
\address{Department of Mathematics \\
         Georgia Institute of Technology \\
         Atlanta, GA 30332-0160, USA}
\email{letu@math.gatech.edu}

\thanks{The author is supported in part by National Science Foundation
 \\
\newline
1991 {\em Mathematics Classification.} Primary 57N10. Secondary 57M25.
\newline
}

 \date{\today \hspace{0.5cm} First edition: November 21, 2005.}


\begin{abstract} We prove that the quantum $SO(3)$-invariant of an
arbitrary 3-manifold $M$ is always an algebraic integer, if the
order of the quantum parameter is co-prime with the order of the
torsion part of $H_1(M,\BZ)$. An even stronger integrality, known as
cyclotomic integrality, was established by Habiro for {\em integral}
homology 3-spheres.  Here we also generalize Habiro's result to all
{\em rational} homology 3-spheres.
\end{abstract}

\maketitle


\addtocounter{section}{-1}


\section{Introduction}
\label{sec.intro}

\subsection{Integrality at roots of non-prime order} Let $\tau_M(q)$
be the quantum $SO(3)$-invariant of a 3-manifold $M$, which can be
defined when $q$ is a root of 1 of {\em odd} order. The quantum
$SU(2)$-invariant was defined by Reshetikhin and Turaev (see
\cite{Turaev}) and the $SO(3)$-version was defined by Kirby-Melvin
\cite{KM} and Turaev. One important result in quantum topology,
first proved by H. Murakami \cite{Murakami}, is that $\tau_M(q)\in
\BZ[q]$, when the order of $q$ is an odd prime and $M$ a rational
homology 3-sphere.
 The result was generalized to all
3-manifolds by Masbaum and Roberts \cite{MR}, using a short,
beautiful proof. Masbaum and Wenzl \cite{MW}, and Takata and Yokota
\cite{KT} generalized the result to some other Lie algebras,
including the $sl_n$ series. The author eventually gave a unified
proof of the integrality for all Lie algebras \cite{Le2}.

The integrality has many important applications, among them is the
construction of of an integral topological quantum field theory of
P. Gilmer and G. Masbaum \cite{Gilmer} and representations of
mapping class groups over $\BZ$.

In the above integrality, the order of $q$ must be an {\em odd
prime}, while  $\tau_M(q)$ can be defined at any root of odd order.
Using a quite different method, Habiro \cite{Habiro1} showed the
integrality of $\tau_M$ when $M$ is an {\em integral homology
3-sphere} at roots of {\em any order}. Thus the restriction on the
order of $q$ is removed, but there is a restriction on the manifold:
$M$ is an integral homology 3-sphere. Then Habiro and the author
proved the integrality of $\tau_M$ for integral homology 3-spheres
at any root of 1  for {\em all simple Lie algebra} \cite{HL}.

The first main result of this paper is to establish the integrality
of $\tau_M$ for {\em all} 3-manifolds with a very minor restriction
on the order of the roots of unity.

\begin{theorem} (a) Suppose $\xi$ is a root of unity of order odd and co-prime
with the order of the torsion part of $H_1(M,\BZ)$, where $M$ is a
closed oriented 3-manifold. Then $\tau_M(\xi) \in \BZ[\xi]$.

(b) Suppose, in addition, $M$ has 0 first Betti number, and $L$ is a
colored, algebraically split link in $M$. Then the quantum invariant
$\tau_{M,L}(\xi) \in \BZ[\xi]$. \label{main1}
\end{theorem}

We will recall the definition of $\tau_M, \tau_{M,L}$ in section
\ref{defs}. In the proof of the theorem we will make use of a result
of G. Andrews concerning Bailey pairs in the theory of $q$-series.

\subsection{Stronger integrality: cyclocomic integrality}

We use the standard notation $(a;q)_n$ of $q$-calculus

$$(a)_n = (a;q)_n := \prod_{i=0}^{n-1}(1-aq^i).$$
 For example, $(q;q)_n = (1-q)(1-q^2)\dots (1-q^n)$. Define the
\newcommand{\Ha}{\widehat{\Lambda}}
Habiro ring $\Ha$ by

$$\Ha := \lim_{\leftarrow n}\Z[q^{\pm 1}]/((q;q)_n).$$

Formally, $\Ha$ is the set of all series of the form

$$f = \sum_{n=0}^\infty f_n(q) \, (1-q)(1-q^2)\dots(1-q^n), \qquad \text{where } \quad f_n(q) \in \Z[q^{\pm 1}].$$

 When $q=\xi$, a root of unity, only a finite number of terms in the right hand side are not 0,
 hence the right hand side defines a complex number, denoted by
 $\ev_\xi(f)$. It is clear that $\ev_\xi(f)\in \BZ[\xi]$.
Thus one can consider every $f\in \Ha$ as a function with domain
 the set of roots of unity. It turns out $\Ha$ has remarkable
properties and plays an important role in quantum topology: Habiro
\cite{Habiro2}  showed that  every $f\in \Ha$ has a Taylor expansion
$T_1(f) \in \BZ[[q-1]]$ which uniquely determines $f$. Also, if
$f(\xi)=g(\xi)$ at infinitely many roots $\xi$ of prime power
orders, then $f=g$ in $\Ha$. The above properties suggest to
consider $\Ha$ as a class of ``analytic functions" with domain the
set of roots of unity.

\begin{theorem} {\rm  (Habiro \cite{Habiro1}) }
 For every integral homology
3-sphere $M$ there is an invariant $I_M\in \Ha$ such that if $\xi$
is a root of unity  then $\ev_\xi(I_M)= \tau_M(\xi)$. \label{000}
\end{theorem}

The integrality of the quantum invariant follows immediate from the
theorem, since for every $f\in \Ha$ and $\xi$ a root of unity,
$f(\xi)\in \BZ[\xi]$. However, belonging to $\Ha$ is a much stronger
integrality. We call it the cyclotomic integrality. For example, the
cyclotomic integrality shows that quantum invariants at infinitely
many roots of unity of prime orders determine the values at any
other roots of unity.  In joint work with Habiro \cite{HL}, we
generalize Theorem \ref{000} to the case of all simple Lie algebras
(but still for integral homology 3-spheres).

\subsection{The case of rational homology 3-spheres} One main result of this paper is to establish
cyclotomic integrality for rational homology 3-spheres.

For a positive integer $d$  let $ A_d := \BZ[\frac{1}{d}][q^{\pm
1/d}]$ and $\BN_d$ the set of positive integers co-prime with $d$.
Denote by $\Phi_s(t)$ the $s$-th cyclotomic polynomial. Let $\Gd
\subset \BQ(q^{1/d})$ be the ring obtained from $A_d$ by adding the
inverses of each $\Phi_s(q^{1/d})$, with $s$ {\em not co-prime with}
$d$:

$$\Gd := A_d[\frac{1}{\Phi_{s}(q^{1/d})}, s \not \in \BN_d].$$
 The ring that replaces Habiro's ring in the case of rational
 homology 3-spheres is

\newcommand{\hRd}{\widehat {\Gd}}
$$\hGd := \lim_{\leftarrow n}\Gd/((q;q)_n)$$

Let $U_d$  be the set of all complex roots of unity with orders {\em
odd} and {co-prime} with $d$. The ring $\hGd$ will play the role of
the Habiro ring, with the set of all roots of 1 replaced by $U_d$.

We first define, for each $\xi\in U_d$, the evaluation map
$\ev_\xi$, which replaces $q$ by $\xi$. Suppose $f\in \BQ[q^{\pm
1/h}]$, where $h$ is co-prime with $r$, the order of $\xi$. There
exists an integer $b$, unique modulo $r$, such that
$(\xi^b)^{h}=\xi$. Then we define

$$\ev_\xi f := f|_{q^{1/h}= \xi^b}.$$

The definition extends to $\ev_\xi: \hGd \to \BC$, since
$\ev_\xi((q;q)_n)=0$ if $n \ge r$. The following is a generalization
of Habiro's result.

\begin{theorem} Suppose $M$ is a rational homology 3-sphere with
$|H_1(M,\BZ)|=d$, a positive integer. There is an invariant $I_M\in
\hGd$ such that if $\xi \in \Ud$, then $\left( \frac{d}{r}\right)
\ev_\xi(q^{(1-d)/4}I_M) $ is the quantum $SO(3)$-invariant of $M$.
Here $\left( \frac{d}{r}\right)$ is the Jacobi symbol. \label{main2}
\end{theorem}

\begin{remark} The reason we single out the factor $q^{(d-1)/4}$ is because it might not be in
$A_d=\BZ[1/d][q^{\pm 1/d}]$; note that it is always in $\BZ[q^{\pm
1/2d}]$. Habiro observed that our proof actually showed that $I_M
\in \hat\Lambda_{d'}$, where $d'$ is the maximal order of elements
of $H_1(M,\BZ)$.
\end{remark}
The cyclotomic integrality has many rigidity properties described in
the next subsections.

\subsection{Evaluation and Taylor expansion of $\hGd$}

For a subset  $\Omega$ of $\Ud$ let

$$ \ev_\Omega: \hGd \to \prod_{\xi\in \Omega} \BC, \qquad \text{defined
by}\ \ev_\Omega(f) = (\ev_\xi(f))_{ \xi\in \Omega}.$$

\begin{theorem} Suppose $\Omega\subset U_d$ contains infinitely many  roots
of order prime powers. If $f,g\in \hGd$ such that
$\ev_\xi(f)=\ev_\xi(g)$ for every $\xi\in \Omega$, then $f=g$. In
other words, the map $\ev_\Omega$ is injective. \label{main3}
\end{theorem}
 We will prove that for every $f\in \hGd$ and $\xi\in \Ud$, one has
$\ev_\xi(f)\in \BZ[1/d,\xi]$. Hence the image of $\ev_\Omega$ is in
$\prod_{\xi\in \Omega} \BZ[1/d,\xi]$. If $f=I_M$ for a rational
homology 3-sphere, then Theorem \ref{main1} shows that $\ev_\xi(f)
\in \BZ[\xi]$.

Any element $f\in \Gd$, considered as a function of the variable
$q^{1/d}$, is analytic at $q^{1/d}=1$. The Taylor series of $f$,
which is a formal power series in $(q^{1/d}-1)$, can be converted
into a formal power series in $(q-1)$ by

$$ q^{1/d}-1 = (1+(q-1))^{1/d} -1 = \sum_{n=1}^\infty
\binom{1/d}{n} (q-1)^n.$$

Thus we obtained an algebra homomorphism $ T_1 : \Gd \to
\BC[[q-1]]$, which can be easily extended to $ T_1 : \hGd \to
\BC[[q-1]]$. We call $T_1(f)$ the Taylor expansion of $f$ at $q=1$,
although in general  there is no analytic continuation of $f$ to a
neighborhood of 1.

\begin{theorem} (a) For any $f\in \hGd$, $T_1(f)$ has
coefficients in $\BZ[1/d]$, i.e. $T_1(f) \in \BZ[1/d][[q-1]]$.

(b) The Taylor expansion map at 1, $T_1: \hGd \to \BZ[1/d][[q-1]]$
is injective. In other words, an element in $\hGd$ is uniquely
determined by its Taylor series.

(c) For a rational homology 3-sphere $M$, the Taylor expansion
$T_1(q^{(1-d)/4}I_M)$ is equal to the Ohtsuki series of $M$
\cite{Ohtsuki1}. \label{main4}
\end{theorem}

Note that because of the factor $q^{(1-d)/4}$, the coefficients of
the Ohtsuki series in general are in $\BZ[1/2d]$. There is also a
similar result for Taylor expansion at any root $\xi\in \Ud$.

The above properties suggests to consider $\hGd$ as a class of
``analytic functions" with domain $\Ud$.

\begin{remark} Actually, we will prove, as suggested by Habiro, that $\hGd$ is
isomorphic to a ring already defined by Habiro, defined without any
denominators. Then Theorems \ref{main3} and \ref{main4} part (b)
follow easily from \cite{Habiro2}. In the original version of the
paper, we had a longer proof, which was also adapted from
\cite{Habiro2}.
\end{remark}

\subsubsection{Applications}

\begin{corollary} Suppose $M$ is a rational homology 3-sphere with
$|H_1(M,\BZ)|=d$.

(a)  The Le-Murakami-Ohtsuki (LMO) invariant of $M$ (see
\cite{LMO,Le3,Ohtsukibook}) determines the quantum $SO(3)$-invariant
at any root $\xi \in \Ud$.

(b) The Ohtsuki series of $M$ has coefficients in $\BZ[1/2d]$.

 (c) If the order of $\xi\in \Ud$ is $p^e$, where $p$ is an odd
prime, then the Taylor series $T_1(q^{(1-d)/4}I_M)$, with $q=\xi$,
converges $p$-adically to $\left (\frac{d}{r}\right) \tau_M(\xi)$.

(d) The values of $\tau_M$ at an infinite subset of roots of order
prime powers in $U_d$ determine the value of $\tau_M$ at any other
root in $\Ud$.
\end{corollary}

Note that part (b) and (c), with $e=1$, were conjectured by R.
Lawrence and first proved by Rozansky \cite{Rozansky} using a quite
different method. Here they, together with part (d), are easy
consequences of our main results. Part (a), which demonstrates one
more universal property \cite{Le3} of the LMO invariant, follows
from the fact that the Ohtsuki series is obtained from the LMO
invariant by way of the $sl_2$-weight system, see
\cite{Ohtsukibook}.

\subsection{Plan of the paper} Section \ref{defs} is devoted to the
definition of quantum invariants, mainly to fix the normalization.
In section \ref{300} we review basic facts of Habiro's work on the
cyclotomic expansion of the colored Jones polynomial. Proofs of
Theorems \ref{main2}, \ref{main3}, \ref{main4}, and their
generalizations are given in sections \ref{401}, \ref{400}, and
\ref{500}.

\subsection{Acknowledgement}
The author wishes to thank A. Beliakova, S. Garoufalidis, T. Ohtsuki
and especially K. Habiro for valuable comments. A great deal of this
paper is inspired by Habiro's work.

\section{Quantum invariants}

\label{defs}

We will use the following notations for elements in $\BZ[q^{\pm
1/2}]$.

$$ \{n\} = q^{n/2}-q^{-n/2},
 \quad  \{n\}!=
\prod_{i=1}^n \{i\} ,\quad  [n] =\frac{\{n\}}{\{1\}}, \quad
\qbinom{n}{k} = \frac{\{n\}!}{\{k\}!\{n-k\}!}.$$

All 3-manifolds in this paper are supposed to be closed and
oriented. Every link in a 3-manifold is framed, oriented, and has
components ordered.
\subsection{The colored Jones polynomial}

\newcommand{\RR} {\mathbf R}

 Suppose $L$ is framed, oriented link
in $S^3$ with $m$ ordered components. For finite-dimensional
$sl_2$-modules $W_1,\dots,W_m$, one can define the quantum invariant
$J_L(W_1,\dots,W_m) \in \BZ[q^{\pm 1/4}]$, see
\cite{Turaev,KM,Ohtsukibook}. The modules $W_1,\dots,W_m$ are
usually called the colors of the link $L$.

It is known that for every positive integer $n$ there is a unique
irreducible $sl_2$-module $V_n$ of dimension $n$. For positive
integers $n_1,\dots,n_m$ we define $ J_L(n_1,\dots,n_m) :=
J_L(V_{n_1},\dots,V_{n_m})$.

We recall here a few well-known formulas, at the same times fix our
normalization. For the unknot $U$ with 0 framing one has

\begin{equation} J_U(n) = [n]= \{n\}/\{1\}. \label{unknot}
\end{equation}

 If $L'$ is obtained from $L$
by increasing the framing of the $i$-th component by 1, then

\begin{equation}\label{framing}
J_{L'}(n_1,\dots,n_m) = q^{(n_i^2-1)/4} J_{L}(n_1,\dots,n_m).
\end{equation}

When all the $n_i$'s are equal to 2, then $J_L(n_1,\dots,n_m)=\tilde
J_L\in \BZ[q^{1/4}]$, a version of the Jones polynomial
\cite{Jones,Ohtsukibook}, which satisfies the skein relation:

\begin{equation*} q^{1/4} \tilde J_{L^+} - q^{-1/4} \tilde J_{L^-} = (q^{1/2}
-q^{-1/2}) \tilde J_{L^0},
\end{equation*}
where $L^+, L^-, L^0$ have blackboard link diagram identical
everywhere, except for a small ball in which $L_+$ has a positive
crossing, $L_-$ a negative crossing, and $L^0$ a resolution of the
crossing.

In general, $J_{L}(n_1,\dots,n_m) \in \BZ[q^{\pm 1/4}]$. However,
there is a number $a\in \{0,\frac{1}{4},\frac{1}{2},\frac{3}{4}\}$
such that $J_{L}(n_1,\dots,n_m) \in q^a\BZ[q^{\pm 1}]$. For a
formula for $a$ and more general results (for all simple Lie
algebra), see \cite{Le_Duke}. A special case is the following.

\begin{lemma}
Suppose  $L$ has 0 linking matrix,  then  $J_{L}(n_1,\dots,n_m) \in
q^{(n_1+\dots n_m-m)/2}\,\BZ[q^{\pm 1}]$. \label{0001}
\end{lemma}

Let $\RR$ be the Grothendieck ring of finite-dimensional
$sl_2$-modules, tensored by $\BC(q^{1/4})$. As a vector space over
$\BC(q^{1/4})$, $\RR$ is freely spanned by $V_1,V_2,\dots$, but as
an algebra, it is isomorphic to $\BC(q^{1/4})[V_2]$. Using linearity
we can define, for $W_i\in \RR$, $J_L(W_1,\dots,W_m) \in
\BC(q^{1/4})$. Thus $J_L$ is considered as a multi-linear function
from $\RR^m$ to $\BC(q^{1/4})$.

\subsection{Definition of $SO(3)$-invariant of 3-manifolds} Suppose $\xi$ is a root
of unity of {\em odd} order $r$ and $f(q;n_1,\dots,n_m)$ a function
of variables $q$ and integers $n_1,\dots,n_m$. Let

$$ {\sum_{n_i}}^\xi f := \sum_{n_i} \ev_\xi (f),$$
where in the sum  all the $n_i$ run the set of {\em odd} numbers
between $0$ and $2r$. Let

$$ F_L(\xi):= {\sum_{n_i}}^\xi
J_L(n_1,\dots,n_m)\prod_{i=1}^m [n_i].$$

Let $U^\pm$ be the unknot with framing $\pm 1$. It is known that
$F_{U^\pm}(\xi) \neq 0$, see Lemma \ref{un1}.

 Suppose $M$ is an
oriented 3-manifold obtained from $S^3$ by surgery along a framed,
oriented link $L$. (Note that $M$ does not depend on the orientation
of $L$). Let $\sigma_+ $ (respectively, $\sigma_-$) be the number of
positive (resp. negative) eigenvalues of the linking matrix of $L$.
Suppose $\xi$ is a root of unity of odd order $r$. Then the quantum
$SO(3)$-invariant is defined by

\begin{equation*}
\tau_M(\xi) = \tau_M^{SO(3)}(\xi) :=
\frac{F_L(\xi)}{(F_{U^{+}}(\xi))^{\sigma_+}\,
(F_{U^{-}}(\xi))^{\sigma_-} }.
\end{equation*}

Suppose in addition to $L$, in $S^3$ there is another framed link
$L'$ with $s$ components, and surgery along $L$ transforms
$(S^3,L')$ to $(M,L'')$. Then the quantum $SO(3)$-invariant of
$(M,L'')$ is defined by

\begin{equation*}
\tau_{M,L''}(k_1,\dots,k_s;\xi) := \frac{{\Sigma_{n_i}^\xi}
(J_{L\cup L'}(n_1,\dots,n_m,k_1,\dots,k_s)\prod_{i=1}^m
[n_i]}{(F_{U^{+}}(\xi))^{\sigma_+}\, (F_{U^{-}}(\xi))^{\sigma_-} }.
\end{equation*}

For connected sum, one has $ \tau_{M\#N}(\xi) =\tau_{M}(\xi)
\tau_N(\xi). $

\subsection{Gauss sum, Laplace transform, and the value of $F_{U^\pm}(\xi)$}

Recall that $\xi\in \Ud$ has odd order $r$. A variation
$\gamma_d(\xi)$ of the Gauss sum is defined by

$$ \gamma_d(\xi):= {\sum_{n}}^\xi q^{d\frac{n^2-1}{4}}.$$

It is known that  $|\gamma_d(\xi)|= \sqrt r$, and hence  is never 0.

Let $\cL_{d;n}: \BZ[q^{\pm n},q^{\pm 1}] \to \BZ[q^{\pm 1/d}]$ be
the $\BZ[q^{\pm 1}]$-linear operator, called the Laplace transform,
defined by

$$\cL_{d;n}(q^{na}) := q^{-a^2/d}.$$

\begin{lemma} Suppose $\xi\in \Sd$ is a root of 1 of order
$r$ and $f \in \BZ[q^{\pm n},q^{\pm 1}]$. Then

$$ {\sum_{n}}^\xi q^{d\frac{n^2-1}{4}}  f = \gamma_d(\xi) \,
\ev_\xi(\cL_{d;n}(f)).
$$
\label{1000}
\end{lemma}

\begin{proof} It's enough to consider the case when $f= q^{na}$,
with $a$ an integer. This case is proved simply by using the
standard completing the square method, see for example, \cite{Le2}.
\end{proof}

The point is that $\cL_{d;n}(f)$, unlike the left hand side
${\sum_{n}^\xi}  q^{d\frac{n^2-1}{4}} f$, does not depend on $\xi\in
\Ud$, and will help us to define a ``universal invariant". Applying
Lemma \ref{1000} with $d=\pm1$ and $f= [n]^2$, using the unknot
formula \eqref{unknot}, we get the following.

\begin{lemma} For the unknot $U^\pm$ with framing $\pm 1$, one has
$F_{U^{\pm}}(\xi) \neq 0$. Moreover,
$$
F_{U^{\pm}}(\xi) = \mp 2  \gamma_{\pm 1}(\xi)\, \ev_\xi\left( \frac{
q^{\mp 1/2}}{\{1\}}\right).
$$
\label{un1}
\end{lemma}

\subsection{Lens spaces} Let $\left( \frac{d}{r}\right)$ be the Jacobi symbol and $s(d,a)$ the
Dedekind sum. Recall that

$$ s(d,a) := \sum_{i=1}^{|d|-1}\left(\left( \frac{i}{d}
\right)\right) \left(\left( \frac{ia}{d} \right)\right), \quad
\text{where $\left(\left( x \right)\right):= x-\lfloor x \rfloor
-1/2$}.
$$

For co-prime integers $a,d$ with $ d>0$, the $SO(3)$-invariant of
the lens space $L(d,a)$, which is obtained by surgery along the
unknot with {\em rational} framing $d/a$, is
    given by (see \cite{Li})

\begin{equation} \tau_{L(d,a)}(\xi) = \left(
\frac{d}{r}\right)\,\ev_\xi\left( q^{-3 s(d,a)} \frac{ q^{1/2d} -
q^{-1/2d} }{q^{1/2} - q^{-1/2}}\right).\label{lens}
\end{equation}

In particular, $\tau_{L(d,a)}(\xi)$ is invertible in $\BZ[\xi]$.

It is also well-known that every non-zero integer $d$, if $\sn(d)$
is the sign of $d$, then
\begin{equation}
\frac{\gamma_d(\xi)}{\gamma_{\sn(d)}(\xi)} = \left
(\frac{|d|}{r}\right )\,\ev_\xi(q^{(\sn(d)-d)/4}). \label{1055}
\end{equation}

\section{Habiro's cyclotomic expansion of the colored Jones
polynomial}

\label{300}
\subsection{The basis $P'_n$} Recall that $\RR$ is isomorphic to $\BC(q^{1/4})[V_2]$, with $V_1,V_2,\dots$ as a
basis over $\BC(q^{1/4})$. Habiro \cite{Habiro1} defined a new basis
$P'_k, k=0,1,2,\dots$, where

$$ P_k':=
\frac{1}{\{k\}!} \, \prod_{i=1}^{k}(V_2-q^{(2i-1)/2}-q^{-(2i-1)/2})
.
$$

The change from the basis $V_n$ to $P'_k$ is given by

\begin{equation*}
V_n= \sum_{k=0}^{n-1}  \qbinom{n+k}{2k+1}\{k\}!\, P'_{k},
\end{equation*}
where the sum over $k$ can be made from 0 to infinity, since
$\qbinom{n+k}{2k+1}=0$ if $k \ge n$.

 For any  link $L$, using the linearity of $J_L$, one has

\begin{equation}
J_{L}(n_1,\dots,n_m) = \sum_{0\le k_i \le n_i-1} J_L(P'_{k_1}, \dots
, P'_{k_m}) \prod_{i=1}^m \qbinom{n_i+k_i}{2k_i+1} \{k_i\}!
\label{151}
\end{equation}

Since there is a denominator in the definition of $P_k'$, one might
expect that $J_L(P'_{k_1}, \dots , P'_{k_m})$ also has non-trivial
denominator. A difficult and important integrality result of Habiro
\cite{Habiro1}  is

\begin{theorem}\cite[Thm.3.3]{Habiro1}
\label{Hab} If $L$ is algebraically split and zero framed link in
$S^3$, then
$$
J_L(P'_{k_1}, \dots , P'_{k_m}) \in \frac{\{2k+1\}!}{\{k\}!\{1\}}
\,\,\BZ[q^{\pm1/2}] = \qbinom{2k+1}{k} (q^2)_k \,\BZ[q^{\pm1/2}],
$$
where $k=\max\{k_1,\dots, k_m\}$.
\end{theorem}

Thus, $J_L(P'_{k_1}, \dots , P'_{k_m})$ is not only integral, but
also divisible by $\qbinom{2k+1}{k}(q^2)_k $, which, in turn, is
divisible by $(q)_k$.

\subsection{Evaluation at $\xi$}

\begin{lemma} Suppose $\xi$ is a root of unity whose order $r$ is
odd and greater than 2. If $k > (r-3)/2$ then

$$ \ev_\xi\left( \frac{\{2k+1\}!}{\{k\}!\{1\}}\right) = 0.$$
\label{54}
\end{lemma}
\begin{proof} First assume that $ k \ge r-1$.
Note that

$$  \frac{\{2k+1\}!}{\{k\}!\{1\}} = \qbinom{2k+1}{k} \{2\}\{3\}\dots
\{k+1\}.$$

Each factor of the right hand side is polynomial in $q^{1/2}$, and
$r$ is among $\{2,3,\dots, k+1\}$. Since $\ev_\xi\{r\}=0$, the
evaluation $\ev_\xi$ of the left hand side is 0.

It remains to consider the case $k+1 <r \le 2k+1$. Note that

$$  \frac{\{2k+1\}!}{\{k\}!\{1\}}= \frac{\{k+1\} \{k+2\} \dots \{2k+1\}}{\{1\}}.$$

When evaluating using $\ev_\xi$, the denominator is not 0, but the
numerator is, since $r$ is among $\{k+1,k+2,\dots,2k+1\}$.
\end{proof}

From the above lemma and  Theorem \ref{Hab} we have the following.
\begin{corollary} Suppose $\xi$ is a root of unity of order $r$ and
$L$ an algebraically split link with 0-framing on each component.
Then

\begin{equation*}
\ev_\xi (J_{L}(n_1,\dots,n_m))  = \ev_\xi \left(
\sum_{k_1,\dots,k_m=0}^{(r-3)/2} J_L(P'_{k_1}, \dots , P'_{k_m})
\prod_{i=1}^m \qbinom{n_i+k_i}{2k_i+1} \{k_i\}!\right)
\end{equation*}
\label{eq.JL}
\end{corollary}

\section{Integrality of quantum invariants}

\label{401}

\subsection{Technical results}
\newcommand{\ty}{\tilde Y}

\subsubsection{Divisibility of the Laplace transform images}

\begin{proposition} \label{222}
 Suppose $\xi\in \Ud$ is a root of order odd $r$ and $k \le
 (r-3)/2$. Choose an integer $b$ such that $db \equiv 1 \pmod r$.
 Then

\begin{align}
{\sum_{n}}^\xi   \qbinom{n+k}{2k+1} \{k\}! \{n\} & =  2
\ev_\xi\left( q^{(k+1)(k+2)/4} \, (q^{k+2})_{r-k-2}
\right) \notag\\
 {\sum_{n}}^\xi  q^{d\frac{n^2-1}{4}} \,
\qbinom{n+k}{2k+1} \{k\}! \{n\} & = - 2\sn(b) \gamma_d(\xi) \,
\ev_\xi(H(k,-b)), \label{121}
\end{align}
where $H(k,-b)$ is in $\BZ[q^{\pm 1/2}]$.
\end{proposition}

The proof will be given later in this section. The upshot here is
the right hand side of \eqref{121} is divisible by $\gamma_d(\xi)$.
An explicit formula for $ H(k,b)$, defined for any pairs $(k,b)$ of
non-zero integers and not depending on $\xi$, as well as the proof
of the proposition, is given in section \ref{pro}.

\subsubsection{Diagonalizing the linking matrix} A link $L$ in a
3-manifold $M$ is {\em algebraically split} if each component of $L$
bounds an orientable surface which does not intersect any other
component.

\begin{proposition} \label{231}
(a) Suppose $M$ is a 3-manifold  with $|Tor(H_1(M,\BZ))|=d$, with
prime decomposition  $d= p_1^{e_1} \dots p_k^{e_k}$. Then there are
lens spaces $L(d_i,a_i)$ such that each $d_i$ is a power of one of
the $p_j$, and the connected sum of $M$ and these lens spaces can be
obtained by integral surgery along an algebraically split link in
$S^3$.

(b) Suppose in addition $M$ is a rational homology 3-sphere, i.e.
$|H_1(M,\BZ)|=d$. Assume that $L$ is an algebraically split link in
 $M$. Then there is an algebraically split link $\bar L$,
which is the disjoint union of 2 sub links $\bar L_1$ and $\bar
L_2$, in $S^3$ such that surgery along $\bar L_1$ transforms
$(S^3,\bar L_2)$ to $(M',L)$. Here $M'$ is the connected sum of $M$
and several lens spaces $L(d_i,a_i)$ with properties as described in
part (a).
\end{proposition}

The proof will be given later. The proposition, part (a), with $M$ a
rational homology 3-sphere, is a modification of Ohtsuki's lemma
\cite{Ohtsuki1}.

\subsection{The case of diagonal linking matrix}
Suppose $M^3$ is obtained from $S^3$ by surgery along  an
algebraically split $(m+s)$-component link $L$ with integral
framings $d_1,\dots,d_{m+s}$ on the link components, where
$d_{m+1}=\dots d_{m+s}=0$ and all other $d_i$'s are not 0. Let $L^0$
denote the link $L$ with all framings switched to 0. Let
$d=\prod_{i=1}^m |d_i|$. Suppose $\xi\in \Ud $ is a root of unity of
odd order $r$. Choose integer $b_i$ such that $b_id_i>0$ and
$b_id_i=1 \pmod r$, for $i=1,\dots,m$.

\subsubsection{Non-degenerate diagonal linking matrix} \label{51}

First we consider the case when $s=0$, i.e. the linking matrix is
non-degenerate, or $M$ is a rational homology 3-sphere.

By framing formula \eqref{framing},

\begin{align}
J_L(n_1,\dots,n_m)&= J_{L^0}(n_1,\dots,n_m) \prod_{i=1}^m
q^{d_i\frac{n_i^2-1}{4}}, \qquad \text {and hence by  formula \eqref{eq.JL}} \notag\\
\ev_\xi(J_L(n_1,\dots,n_m))&= \ev_\xi \left(\sum_{k_i=0}^{(r-3)/2}
J_{L^0}(P'_{k_1},\dots,P'_{k_1}) \prod_{i=1}^m
q^{d_i\frac{n_i^2-1}{4}}\,\qbinom {n_i+k_i}{2k_i+1} \{ k_i\}!\right)
\notag
\end{align}

By definition,
\begin{align}
F_L(\xi)& =   {\sum_{n_i}}^\xi \left( J_L(n_1,\dots,n_m) \prod_{i=1}^m [n_i] \right), \qquad \text{and hence} \notag\\
 F_L(\xi)& = \ev_\xi\left( \sum_{k_i=0}^{(r-3)/2}
J_{L^0}(P'_{k_1},\dots,P'_{k_m}) \prod_{i=1}^m {\sum_{n_i}}^\xi
q^{d_i\frac{n_i^2-1}{4}} \qbinom {n_i+k_i}{2k_i+1} \frac{\{ k_i\}!
\{n_i\}}{\{1\}}
\right) \label{156} \\
&= \left( \prod_{i=1}^m  \frac{-2\sn(d_i)
\gamma_{d_i}(\xi)}{\{1\}}\right) \, \ev_\xi\left(
\sum_{k_i=0}^{(r-3)/2} J_{L^0}(P'_{k_1},\dots,P'_{k_m})
\prod_{i=1}^m  H(k_i,-b_i) \right) \qquad \text{by \eqref{121}}.
\end{align}

 Using Lemma \ref{un1} and the definition of $\tau_M$, one
obtains

\begin{equation}
\tau_M(\xi) = \prod_{i=1}^m
\frac{\gamma_{d_i}(\xi)}{\gamma_{\sn(d_i)}(\xi)}  \,  \ev_\xi\left(
\sum_{k_i=0}^{(r-3)/2} J_{L^0}(P'_{k_1},\dots,P'_{k_m})
\prod_{i=1}^m  q^{\sn(d_i)1/2} H(k_i,-b_i) \right). \notag
\end{equation}

Using \eqref{1055} for
$\frac{\gamma_{d_i}(\xi)}{\gamma_{\sn(d_i)}}$, we get

\begin{equation}
\tau_M(\xi) =  \left( \frac{d}{r}\right)  \,  \ev_\xi\left(
\sum_{k_i=0}^{(r-3)/2} J_{L^0}(P'_{k_1},\dots,P'_{k_m})
\prod_{i=1}^m  q^{(3\sn(d_i)-d_i)/4} H(k_i,-b_i) \right).
\label{510}
\end{equation}

\subsubsection{Degenerate diagonal linking matrix} Now we assume the
general case, $s$ might not be 0.

The argument that leads to \eqref{156} will give us

\begin{multline*}F_L(\xi) = \sum_{k_i=0}^{(r-3)/2} \ev_\xi\left(
J_{L^0}(P'_{k_1},\dots,P'_{k_{m+s}}) \prod_{i=1}^m {\sum_{n_i}}^\xi
q^{d_i\frac{n_i^2-1}{4}} \qbinom {n_i+k_i}{2k_i+1}
\frac{\{ k_i\}! \{n_i\}}{\{1\}}\right)\times \\
\left(\prod_{i=m+1}^{m+s} {\sum_{n_i}}^\xi \qbinom {n_i+k_i}{2k_i+1}
\frac{\{ k_i\}! \{n_i\}}{\{1\}}\right)  .
\end{multline*}

Using Proposition \ref{222}, we get
\begin{multline}
 \tau_M(\xi) =  \left( \frac{d}{r}\right) \sum_{k_i=0}^{(r-3)/2}  \,
\ev_\xi\left( J_{L^0}(P'_{k_1},\dots,P'_{k_{m+s}}) \prod_{i=1}^m
q^{(3\sn(d_i)-d_i)/4} H(k_i,-b_i) \right) \times \\
\ev_\xi \left (\prod_{i=m+1}^{m+s} 2 q^{(k_i+1)(k_i+2)/4} \, \frac{
(q^{k_i+2})_{r-k_i-2}}{\{1\}}\right)\label{2111}
\end{multline}

Note that, since $k_i\le (r-3)/2$, one has $ \frac{
(q^{k_i+2})_{r-k_i-2}}{\{1\}}\in \BZ[q^{\pm 1/2}]$. It follows that
$\tau_M(\xi)\in \BZ[\xi]$.

\subsection{Proof of Theorem \ref{main1}}
{Part (a).} By the product formula, $\tau_M=\tau_{(M\#N)}/\tau_N$,
if $\tau_N\neq 0$. Since $\tau_M(\xi)$, with $M$ the lens space
$L(d,a)$ and $\xi\in \Ud$, is invertible   in $\BZ[\xi]$, using
Proposition \ref{231} we can assume that $M$ is obtained from a link
with diagonal linking matrix as described in the previous
subsection. There we have proved that $\tau_M(\xi) \in \BZ[\xi]$.

{Part (b).} Again we can assume that $(M,L)$ is obtained from
$(S^3,\bar L_2)$ by surgery along $\bar L_1$, as described in
Proposition \ref{231}.

Similarly to \eqref{2111} we have

\begin{multline*} \tau_{M,L}(n_{m+1}, \dots, n_{m+s}) =  \left(
\frac{d}{r}\right)  \, \ev_\xi\left( \sum_{k_i=0}^{(r-3)/2} J_{\bar
L^0}(P'_{k_1},\dots,P'_{k_{m+s}}) \prod_{i=1}^m
q^{(3\sn(d_i)-d_i)/4} H(k_i,-b_i) \right)\times\\
\ev_\xi\left( \prod_{i=m+1}^{m+s} \qbinom{n_i+k_i}{2k_i+1} \{k_i\}!
\right),
\end{multline*}
which is in $\BZ[\xi]$.\qed

\subsection{Proofs of technical results} \label{pro}

For a non-negative integer $k$ let

$$ Z(k) := \sum_{j=0}^{2k+1}(-1)^j \qbinom{2k+1}{j}\,
t^{(j-k)^2} \quad \in \quad \BZ[q^{\pm 1}, t^{\pm 1}].$$

For integer $d\neq 0$ and arbitrary integer $b$ define the
$\BZ[q^{\pm 1}]$-linear algebra operators

\begin{align*}
\varphi_d: \BZ[q^{\pm 1}, t^{\pm 1}] \to \BZ[q^{\pm 1/d}],  \qquad
& \text{where}\ \varphi_d(t) := q^{1/d},\\
\tilde \varphi_b: \BZ[q^{\pm 1}, t^{\pm 1}] \to \BZ[q^{\pm 1}],
 \qquad
& \text{where}\ \tilde \varphi_b (t) = q^{b}.
\end{align*}

The relation between $\varphi$ and $\tilde \varphi$ is: If $\xi\in
\Ud$ has order $r$, and $db\equiv 1 \pmod r$, then  for every $f\in
\BZ[q^{\pm 1}, t^{\pm 1}]$:

\begin{equation}\ev_\xi(\varphi_d(f)) = \ev_\xi(\tilde \varphi_b(f)).
\label{33}
\end{equation}

For non-zero integer $d$ and arbitrary integer $b$ let us define

\begin{align*}
Y(k,d) & := \varphi_d(Z(k)) =  \sum_{j=0}^{2k+1}(-1)^j
\qbinom{2k+1}{j}\,
q^{(j-k)^2/d}   \\
\ty(k,b) &:= \tilde\varphi_b(Z(k)) =  \sum_{j=0}^{2k+1}(-1)^j
\qbinom{2k+1}{j}\, q^{b(j-k)^2}.
\end{align*}

\begin{lemma} \label{110}
Suppose $bd \equiv 1 \pmod r$, where $r$ is  the order of $\xi\in
\Ud$, then

\begin{align*} {\sum_{n}}^\xi  q^{d\frac{n^2-1}{4}} \,  \qbinom{n+k}{2k+1} \{k\}!
\{n\} & = -2 \gamma_d(\xi) \, \ev_\xi \left(\frac{\tilde Y(k,-b)\{
  k\}!}{\{2k+1\}!}\right).
\end{align*}
\end{lemma}

\begin{proof}  We will first
calculate $\cL_{d;n}(\{n\} \, \{n+k\}!/\{n-k-1\}!)$.

Using $\{j\}= - q^{-j/2} (1-q^j)$, one sees that

\begin{align}
\{n\} \, \{n+k\}!/\{n-k-1\}! &= q^{-nk}(q^{-n} -1) \, (q^{n-k};q)_{2k+1}\notag \\
&= q^{-n-nk} (q^{n-k};q)_{2k+1} - q^{-nk}\, (q^{n-k};q)_{2k+1}.
\label{010}
\end{align}

It is easy to check that the two terms of the right hand side of
\eqref{010} can be obtained from one another by the involution $n\to
-n$. Since $\cL_{d;n}$ is invariant under $n \to -n$, one has
\begin{equation}\cL_{d;n}(\{n\} \, \{n+k\}!/\{n-k-1\}!) = -2 \cL_{d;n}(q^{-nk}\,(q^{n-k};q)_{2k+1}
).\label{123}
\end{equation}

Explicit expansion of $(a;q)_k$ is well-known, and we have

\begin{equation} q^{-nk} (q^{n-k};q)_{2k+1} = \sum_{j=0}^{2k+1}(-1)^j
\qbinom{2k+1}{j} \, q^{n(j-k)}. \label{4}
\end{equation}

Using the definition of $\cL_{d;n}$ and \eqref{123} we get

\begin{align} \cL_{d;n}(\{n\} \, \{n+k\}!/\{n-k-1\}!)
&= -2 Y(k,-d). \notag
\end{align}

Multiplying by $\{k\}!/\{2k+1\}!$ and using Lemma \ref{1000} and
\eqref{33}, we get the lemma.
\end{proof}

\subsubsection{Factoring $\ty(k,b)$} It turns out that $\ty(k,b)$ is always divisible by
$\frac{\{2k+1\}!}{\{k\}!}$. To describe the quotient let us define,
for positive integer $b$,

$$ H(k,b): =
q^{(k+1)(k+2)/4}\sum_{k+1 \ge n_b \ge n_{b-1} \ge \dots \ge n_{2}\ge
0} q^{n_2^2+n_3^3+\dots +n_b^2} \frac{(q)_{k+1}}{\prod_{i=1}^{b}
(q)_{n_{i+1}-n_i}}
 ,$$

where $n_{b+1}=k+1$ and $n_1=0$, and  $(a)_n$ stands for $(a;q)_n$.
For example,

\begin{align} H(k,1) & = q^{3/2}   \notag\\
H(k,2) &= q^3\sum_{j=0}^{k+1} q^{j^2} \frac{(q)_{k+1}}{(q)_{j} (q)_{k+1-j}}  = -q^3\sum_{j=0}^{k+1} q^{ j(j+k+1)/2}\, \qbinom{k+1}{j} \notag\\
H(k,3) &= q^5 \sum_{0\le j\le l \le k+1} q^{j^2 + l^2 }\,
\frac{(q)_{k+1}}{(q)_j \,(q)_{l-j} (q)_{k+1-l}}. \notag
\end{align}

 Let $\overline{H(k,b)}$ be obtained from $H(k,b)$ by the involution $q\to q^{-1}$. For $b=0$ define $H(k,b)=0$ and for  $b<0$ by
$$H(k,-b) := (-1)^k \overline{H(k,b)}.$$

\begin{remark} Habiro observed there is a similarity between
$H(k,b)$ and the coefficient $c'_{n,p}$ of Habiro's twists in
\cite{Masbaum}. We have found the exact relation and will discuss it
in another publication.
\end{remark}

The following theorem is the main technical result.

\begin{theorem}\label{tech} For every  integer
$b$, $\ty(k,b)$ is divisible by $\frac{\{2k+1\}!}{\{k\}!}$. Moreover

$$ \ty(k,b)  \frac{\{k\}!}{\{2k+1\}!} = -\sn(b)  H(k,b).$$

\end{theorem}

\begin{proof} If $b=0$, then $\ty(k,b)$ is the right hand side of
\eqref{4}, with $n=0$. In this case (with $n=0$), the left hand side
of \eqref{4} is 0. Hence  $\ty(k,0)=0$, and we are done.

Suppose $b\neq 0$. Using the involution $q\to q^{-1}$, we can assume
that $b>0$.

Let $\alpha_n,\beta_n$ be a Bailey pair as defined in Section 3.4 of
\cite{Andrews}, with $a=1$. Then for any numbers $b_i,c_i,
i=1,\dots,k$ and positive integer $N$ we have the identity (3.43) of
\cite{Andrews}:

\begin{multline}\sum_{n\ge 0} (-1)^n\alpha_n q^{-\binom{n}{2}+kn+Nn} \frac{(q^{-N})_n}{(q^{N+1})_n} \prod_{i=1}^k
\frac{(b_i)_{n}}{b_i^{n}} \frac{(c_i)_{n}}{c_i^{n}}
\frac{1}{(\frac{q}{b_i})_n (\frac{q}{c_i})_n}=
\\
\frac{(q)_N\,  (\frac{q}{b_kc_k})_N}{(\frac{q}{b_k})_N\,
(\frac{q}{c_k})_N} \sum_{n_k \ge n_{k-1} \ge \dots \ge n_1\ge 0}
\beta_{n_1} \frac{q^{n_k} (q^{-N})_{n_k} (b_k)_{n_k}
(c_k)_{n_k}}{(q^{-N}b_kc_k)_{n_k}}\prod_{i=1}^{k-1}\frac{q^{n_i}
\frac{(b_i)_{n_i}}{b_i^{n_i}} \frac{(c_i)_{n_i}}{c_i^{n_i}}
(\frac{q}{b_i c_i})_{n_{i+1}-n_i} }{(q)_{n_{i+1}-n_i}
(\frac{q}{b_i})_{n_{i+1}}  (\frac{q}{c_i})_{n_{i+1}}}. \label{1002}
\end{multline}

A special Bailey pair is given by (see section 3.5 of
\cite{Andrews}):

\begin{align*} \alpha_0 &=1, \quad \alpha_n = (-1)^n q^{n(n-1)/2}(1+q^n) \quad
\text{for } n \ge 1.\\
\beta_0 & =1, \quad \beta_n =0 \quad \text{for } n \ge 1.
\end{align*}
 Using the obvious limits

\begin{align*} \lim_{c\to \infty} \frac{(c)_n}{c^n} & =(-1)^n q^{n(n-1)/2}\\
\lim_{c\to \infty} \left(\frac{q}{c}\right)_n & =1\\
\lim_{c_1,c_2\to
\infty}\frac{(c_1)_n(c_2)_n}{(q^{-N}c_1c_2)_n}&=(-1)^n q^{n(n-1)/2}
q^{Nn}
\end{align*}
we see that the limit of the left hand side of \eqref{1002} , when
 $b_i,c_i$ tend to infinity, with $k$ replaced by
$b$, is

$$ LHS = 1 + \sum_{n=1}^N \frac{q^{Nn} (q^{-N})_n}{(q^{N+1})_n}\,
(1+q^n) \, q^{bn^2}.$$

Here the first term corresponds to $n=0$, and the sum of the second
term terminates at $n=N$ since $(q^{-N})_n=0$ if $n > N$. It is easy
to check that

$$ \frac{q^{Nn} (q^{-N})_n}{(q^{N+1})_n}\,
(1+q^n) = \frac{(-1)^n}{\qbinom{2N}{N}} \qbinom{2N}{N-n} (q^{n/2} +
q^{-n/2}).$$

Hence, with $N=k+1$, we have

\begin{equation} LHS = 1 + \frac{1}{\qbinom{2k+2}{k+1}}
\sum_{n=1}^{k+1} (-1)^n\qbinom{2k+2}{k+1-n} (q^{n/2} + q^{-n/2})\,
q^{bn^2}. \label{1001}
\end{equation}

Recall that

$$ \ty(k,b) = \sum_{j=0}^{2k+1}(-1)^j \qbinom{2k+1}{j}\, q^{b
(j-k)^2}  .$$

The sum can be assume from $j=-1$ to $j=2k+1$, since the term of
$j=-1$ is equal to 0. Separating the case $j=k$ and combining any
other $j$ with $2k-j$, we have

$$ \ty(k,b) = (-1)^k  \qbinom{2k+1}{k} + \sum_{j=-1}^{k-1}(-1)^jq^{b
(j-k)^2} \left(  \qbinom{2k+1}{j} + \qbinom{2k+1}{2k-j} \right) .$$

It is easy to check that

$$ \qbinom{2k+1}{j} + \qbinom{2k+1}{2k-j} =
\frac{\{k+1\}}{\{2k+2\}}\qbinom{2k+2}{j+1}(q^{(k-j)/2} +
q^{(j-k)/2}) .$$

Hence, using new parameter $n= k-j$,  we see that

$$ \ty(k,b) = (-1)^k  \qbinom{2k+1}{k} +  \frac{\{k+1\}}{\{2k+2\}} \sum_{n=1}^{k+1}(-1)^{k-n}q^{b
n^2} \qbinom{2k+2}{k+1-n} (q^{n/2} + q^{-n/2}) .$$

Using $ \qbinom{2k+2}{k+1} =
\qbinom{2k+1}{k}\frac{\{2k+2\}}{\{k+1\}}$, from \eqref{1001} we get,

\begin{equation}
\ty(k,b) = (-1)^{k} \qbinom{2k+1}{k} LHS. \label{001}
\end{equation}

The limit of the right hand side of \eqref{1002}, when $b_i,c_i$
tend to infinity, with $k$ replaced by $b$, is

\begin{equation*}  RHS = (q)_N \sum_{n_b\ge n_{b-1}\ge \dots \ge
n_1=0}(-1)^{n_b} q^{n_b(n_b-1)/2+ N n_b+n_b} (q^{-N})_{n_b}
\prod_{i=1}^{b-1} \frac{q^{n_i^2}}{(q)_{n_{i+1}-n_i}}
\end{equation*}
Note that $n_1=0$, since this is the only case when $\beta_{n_1}\neq
0$, and in the sum the index $n_b\le N$, since $(q^{-N})_{n_b}=0$ if
$n_b >N$. An easy calculation shows that, with $k=N-1$,

\begin{equation} RHS = (-1)^k\{k+1\}! H(k,b). \label{002}
\end{equation}

Since $LHS=RHS$, from \eqref{001} and \eqref{002} we get $
\displaystyle{\ty(k,b) = \frac{\{2k+1\}!}{\{k\}!}H(k,b)}$.
\end{proof}

\subsubsection{Proof of Proposition \ref{222}}

Part (b), the difficult part, follows from Lemma \ref{110} and
Theorem \ref{tech}.

Now we prove part (a). Again noting  that the two terms of the right
hand side of \eqref{010} can be obtained from one another by the
involution $n\to -n$, and using \eqref{4}, we have

\begin{equation*}{\sum_{n}}^\xi  \{n\} \{n+k\}!/\{n-k-1\}! = -2
{\sum_{n}}^\xi  \left( \sum_{j=0}^{2k+1}(-1)^j \qbinom{2k+1}{j} \,
q^{n(j-k)} \right) .
\end{equation*}

Note that ${\sum_{n}^\xi}  (q^{na})= r$ or $0$, according as $a$ is
divisible by $r$ or not. With $0\le j \le 2k+1 \le r-2$, the only
case with $j-k$ divisible by $r$ is when $j=k$. Hence

$$ {\sum_{n}}^\xi   \{n\} \{n+k\}!/\{n-k-1\}!= -2(-1)^k r \times
\ev_\xi \qbinom{2k+1}{k}.$$

Multiplying both side by $\frac{\{k\}!}{\{2k+1\}!}$ and using the
well-known $r= \prod_{i=1}^{r-1}(1-\xi^i)$, we get part (a). \qed

\subsection{Proof of Proposition \ref{231}} Here we modify
Ohtsuki's proof of a similar result \cite{Ohtsuki1}.

\subsubsection{Linking pairing}
A {\em linking pairing} on a finite abelian group $G$ is a
non-singular symmetric bilinear map from $G\times G$ to $\BQ/\BZ$.
Two linking pairing $\nu, \nu'$ on respectively $G,G'$ are
isomorphic if there is an isomorphism between $G$ and $G'$ carrying
$\nu$ to $\nu'$. With the obvious block sum, the set of equivalence
classes of linking pairings is a semigroup.

One type of linking pairing is given by  non-singular square
symmetric matrices with integer entries: any such $n\times n$ matrix
$A$  gives rives to a linking pairing $\phi(A)$ on $G= \BZ^n /A
\BZ^n$ defined by $\phi(A)(v,v') = v^t A^{-1} v' \in \BQ \mod \BZ$,
where $v,v'\in \BZ^n$. If there is a {\em diagonal} matrix $A$ such
that a linking pairing $\nu$ is isomorphic to $\phi(A)$, then we say
that $\nu$ is {\em of diagonal type}.

Another type of pairing is the  pairing $\phi_{b,a}$, with $a,b$
non-zero co-prime integers, defined on the cyclic group $\BZ/b$ by
$\phi_{b,a}(x,y)= axy/b \mod \BZ$. It is clear that $\phi_{b,\pm 1}$
is also of the former type, namely, $\phi_{b,\pm 1}= \phi(\pm b)$,
where $(\pm b)$ is considered as the $1\times 1$ matrix with entry
$\pm b$.

\begin{proposition} Suppose $|G|=d$, with prime decomposition $d=\prod_{i=1}^k p_i^{e_i}$, and
$\nu$  a linking pairing on $G$. There are pairs of non-zero,
co-primes integers $(b_j,a_j),j=1,\dots,s$ such that each $b_j$ is a
power of some $p_i$, and the block sum of $\nu$ and all the
$\phi_{b_j,a_j}$ is of diagonal type. \label{1101}
\end{proposition}

\begin{proof} The semigroup of linking pairing has the following
 generators in 3 groups, see \cite{Kojima,Wall}:

Group 1: $\phi(\pm p^k)$, where $p$ is a prime, and $k>0$.

Group 2: $\phi_{b,a}$ with $b=p^k$ as in group 1, and $a$ is a
non-quadratic residue modulo $p$ if $p$ is odd, or $a=\pm 3$ if
$p=2$.

Group 3: $E_0^k$  on the group $\BZ/2^k \oplus \BZ/2^k$ with $k \ge
1$ and $E_1^k$ on the group $\BZ/2^k \oplus \BZ/2^k $ with $k \ge
2$.

For explicit formulas of $E_0^k$ and $E_1^k$, see \cite{Kojima}. We
will use only a few relations between these generators, taken from
\cite{Kojima,Wall}. It's enough to prove the lemma when $\nu$ is one
of the generator.

Any pairing in group 1 is already of the form $\phi(\pm  p^k)$.

Let $\nu=\phi_{b,a}$ be in group 2. Suppose  $p$ is odd, then one of
the relations is $\phi_{b,a}\oplus \phi_{b,a} = \phi(b) \oplus
\phi(b)$, which is diagonal type. Suppose $b=2^k$, then $a=\pm 3$,
and one of the relations says $\phi_{b,\pm 3}\oplus \phi_{b,\pm 3} =
\phi(\mp b) \oplus \phi(\mp b).$

Suppose $\nu=E_0^k$. One of the relations is $E^k_0 \oplus
\phi(-2^k) = \phi(2^k)\oplus \phi(-2^k) \oplus \phi(-2^k)$.

Finally let $\nu= E_1^k$. One of the relations is $E_1^k \oplus
\phi_{2^k,3}= \phi(2^k)\oplus \phi(2^k)\oplus \phi(2^k)$.
\end{proof}

\subsubsection{Linking pairing on the torsion group of $H_1(M,\BZ)$}
Suppose $M$ is obtained from $S^3$ by surgery along
a framed oriented link $L$, with non-degenerate linking matrix $A$.
 Then $M$ is a rational homology
3-sphere, and the linking pairing on $H_1(M,\BZ)$ is exactly
$\phi(A)$. The following proposition was already implicitly given in
\cite{Ohtsuki1}.

\begin{proposition}  If the linking pairing of a rational homology
3-sphere $M$ is of diagonal type, then $M$ can be obtained from
$S^3$ by surgery along an algebraically split link with integer
framings. \label{111}
\end{proposition}

\begin{proof}

Suppose $M$ is obtained from $S^3$ by surgery along a framed link
$L$ with linking matrix $A$. By assumption, $\phi(A) \cong \phi(B)$,
where $B$ is a diagonal matrix. In this case, it is known that there
is a unimodular integral matrix $P$ such that $P^t A'P=B'$, where
$A'$ and $B'$ are obtained from respectively $A$ and $B$ by
block-adding a  diagonal matrix with $\pm 1$ on the diagonal. Using
the two Kirby moves on links, one can easily go from $L$ to another
framed link with linking matrix $B'$, which is diagonal. Surgery on
the new link yields the same 3-manifold.
\end{proof}

\subsubsection{Proof of Proposition \ref{231}}

\newcommand{\Tor}{\mathrm{Tor}}

Part (a). Suppose $H_1(M,\BZ)= \BZ^r \oplus \Tor$, where $\Tor$ is
the torsion part. Choose disjoint curves $\alpha_1,\dots,\alpha_r$
in $M$ representing generators of the infinite part $\BZ^r$ of
$H_1(M,\BZ)$. By Poincare duality, there are oriented surfaces
$S_1,\dots S_r$ in $M$ such that the algebraic intersection number
between $\alpha_i$ and $S_j$ is $\delta_{ij}$. Using the standard
tube construction if needed, we can assume that $S_i$ meets
$\alpha_i$ at exactly 1 point, and misses all other $\alpha_j$.

Let $N(\alpha_i)$ be a small tubular neighborhood of $\alpha_i$ in
$M$, and $S'_i= S_i \setminus N(\alpha_i)$. Then  $\beta_i:=
\partial S_i'$ is the intersection of $S_i$ with $\partial
N(\alpha_i)$. Removing the interior of $N(\alpha_i)$ from $M$ and
regluing back in such a way that the Dehn filling kills the homology
class $\alpha_i$, we get a rational homology 3-sphere $M'$, together
with $S'_i$ in $M'$. Each $S'_i$ bounds the curve $\beta_i$ in $M'$,
such that doing surgery along $\beta_i$ with framing 0, from $M'$ we
get back $M$. Note that $H_1(M',\BZ)=\Tor$.

If $N= L(b,a)$, the lens space, then the linking pairing on
$H_1(N,\BZ)$ is exactly $\phi_{b,a}$. By Proposition \ref{1101},
there are lens spaces $L(d_i,a_i)$, with properties as described in
the statement of Proposition \ref{231} , such that the linking
pairing of $M'' := M'\#( \#_i L_{d_i,a_i})$ is of  diagonal type. By
Proposition \ref{111}, $M''$ can be obtained from $S^3$ by integral
surgery along a link $L' \subset S^3$ with diagonal linking matrix.
Thus $S^3 \setminus N(L') = M'' \setminus N(L'')$, for a link $L''$
in $M''$. One can isotope $L''$ in $M''$ off the surfaces $S_i'$,
since the surface with non-trivial boundary can be isotoped to a lie
in a small neighborhood of a graph. Thus all the surfaces $S_i'$ can
be considered as lying in $S^3\setminus L'$. The link $L'\cup \{
\beta_1,\dots,\beta_r \}$ has diagonal linking matrix, with framing
0 on $\beta_i$. Doing surgery on $L'\cup \{ \beta_1,\dots,\beta_r
\}$ gives us $M'$. This completes the proof of part (a). The proof
of part (b) is totally similar.\qed

\section{Universal invariant in cyclotomic rings} \label{400}
\subsection{ Proof of Theorem \ref{main2}} First notice that if one
can find an element $f\in \hGd$ such that

\begin{equation}
\left( \frac{d}{r}\right) \tau_{M}(\xi) = \ev_\xi ( q^{(1-d)/4}f),
\label{5555}
\end{equation}
then the injectivity of the map $\ev_\Omega$, with $\Omega=\Ud$ in
Theorem \ref{main3}, would show that $f$ is an invariant of $M$. The
task now is to find such an $f$ for every rational homology
3-sphere. Let us consider 3 cases.

Case 1: $M=L(d,a)$, a lens space, with $d>0$. Let

$$ I_{M} :=    q^{3s(d,1)-3 s(d,a)}
\frac{1-q^{-1/d}}{1-q^{-1}}.$$

It is well-known that $3d (s(d,1)- s(d,a)) \in \BZ$, see
\cite{Borwein}, hence $I_M \in \Lambda_d$. Using
$s(d,1)=(d-1)(d-2)/12d$ and \eqref{lens}, one gets \eqref{5555}.
 Moreover, $I_{L(d,a)}$ is invertible in $\Lambda_d$.

Case 2: $M$ is obtained from $S^3$ by integral surgery along an
algebraically split link $L$ as in section \ref{51}. We will use the
notation of section \ref{51}. The following proposition will be
proved later in this section.

\begin{proposition} The element
$ \displaystyle{ \frac{Y(k,d) \, \{k\}!}{\{2k+1\}!}}$  is in
$q^{(k+1)(k+2)/4}\Lambda_d.$ \label{5}
\end{proposition}

Note that, by \eqref{33} and Theorem \ref{tech},

\begin{align}
\sn(b_i) \ev_\xi \left(Y(k_i,-d_i) \frac{\{k_i\}!
}{\{2k_i+1\}!}\right)&= H(k_i,-b_i) . \label{53}
\end{align}

 Let

\begin{equation}
I_M:= q^{(d-1)/4}\sum_{k_i=0}^{\infty}
J_{L^0}(P'_{k_1},\dots,P'_{k_m}) \prod_{i=1}^m \sn(d_i)
q^{(3\sn(d_i)-d_i)/4} Y(k_i,-d_i) \frac{\{k_i\}! }{\{2k_i+1\}!}
\end{equation}

We will prove later the following.
 \begin{lemma} One has $I_M \in \hGd$.
 \label{6}

\end{lemma}

  By Habiro's theorem (Theorem
\ref{Hab}) $J_{L^0}(P'_{k_1},\dots,P'_{k_m})$ is divisible by
$\frac{\{2k+1\}!}{\{k\}!\{1\}}$. Hence by Lemma \ref{54}, if $k_i >
(r-3)/2$ for some $i$, then
$\ev_\xi(J_{L^0}(P'_{k_1},\dots,P'_{k_m}))=0$. Thus combining with
\eqref{53} and \eqref{510}, we get \eqref{5555}.

Case 3: $M$ is an arbitrary rational homology 3-sphere. By the
diagonalizing lemma, there are lens spaces $M_1,\dots,M_l$ such that
$M'= (\#_{i=1}^l M_i) \# M$ is of the case 2. Since each $I_{M_j}$
is invertible in $\Gd$, we can define

$$ I_{M} = I_{M'} ( \prod_{i=1}^l I_{M_i})^{-1}.$$

Using the product formula, we see that $I_M$ satisfies \eqref{5555}.
This completes the proof of Theorem \ref{main2}.

\subsection{Proof of Lemma \ref{6}} By Proposition \ref{5},
$Y(k_i,-d_i) \frac{\{k_i\}! }{\{2k_i+1\}!}$ is in
$q^{(k_i+1)(k_i+2)/4}\Gd$, and by Theorem ~\ref{Hab},
$J_{L^0}(P'_{k_1},\dots,P'_{k_m})$ is divisible by $(q;q)_n$. The
only problem is, a priori, the term in the sum formula of $I_M$
might contain the fractional power $q^{1/2d}$, and all we need to
show is we need only $q^{1/d}$, but not $q^{1/2d}$.

First, using Lemma \ref{0001} one sees that
$\displaystyle{J_{L^0}(P'_{k_1},\dots,P'_{k_m})  \prod_{i=1}^m
\frac{q^{\sn(d_i)/2} \{k_i\}!}{\{2k_i+1\}!}}$  is in $ \BZ(q)$,  (no
fractional power of $q$.) It suffices to prove that $q^{(d-1)/4}
\prod_{i=1}^m q^{(\sn(d_i) - d_i)/4}$ is in $\BZ[q^{\pm 1/d}]$. This
is equivalent to the fact that

\begin{equation*}
D:= d(d-1) -d \sum_{i=1}^m (\sn(d_i) - d_i)) \quad \text{ is
divisible by $4$}.
\end{equation*}

We have  $ \sn(d_i) -d_i = \sn(d_i)(1-|d_i|)$. If the sign
$\sn(d_i)$ is changed, then $D$ is altered by $\pm 2d(1-|d_i|)$,
which is divisible by $4$, since $d_i|d$. Hence we can assume that
every $\sn(d_i)=+1$.

We use induction on $m$. If $m=1$, then the statement is trivial.
Note that

$$(1-d_1d_2) -(1-d_1) - (1-d_2) = (1-d_1)(1-d_2),$$
which, after multiplied by $d_1d_2$, is divisible by 4. Hence the
induction step can be carried out by replacing $d_1$ and $d_2$ with
$d_1+d_2$. \qed

\newcommand{\cI}{\mathcal I}
\subsection{Proof of Proposition \ref{5}} For every non-negative
integer $n$ let $\cI_n$ be the ideal in $\BZ[t^{\pm 1}, q^{\pm1}]$
generated by $c_{k,n}:=(t;q)_k (q^{k+1})_{n-k}, k=0,1,\dots,n$. This
ideal was used by Habiro in \cite{Habiro3}. Recall that
$\tilde\varphi_b: \BZ[t^{\pm 1}, q^{\pm1}] \to \BZ[q^{\pm1}]$ is the
$\BZ[ q^{\pm1}]$-algebra homomorphism defined by
$\tilde\varphi_b(t)= q^b$.

\begin{proposition} Suppose  $f(t,q) \in \BZ[t^{\pm1},q^{\pm1}]$.
Then $\tilde\varphi_b(f)$ is divisible by $(q)_n$ for every integer
$b$ if and only if $f \in \cI_n$.
 \label{333}
\end{proposition}

\begin{proof}
We use induction on $n$. Suppose the statement holds true for $n-1$.
By induction,

$$ f= a_0(t,q) c_{0,n-1} + a_1(t,q) c_{1,n-1} +\dots +a_{n-1}(t,q) c_{n-1,n-1}.$$
Applying $\tilde \varphi_0$, noting that $\tilde
\varphi_0(c_{k,n-1})=0 $ if $k >0$, one gets that $a_0(1,q)$ is
divisible by $1-q^n$, or $a_0(t,q) \in (1-t) + (1-q^n)$. Note that
$(1-t)c_{0,n-1}+(1-q^n)c_{0,n-1} \subset \cI_n$. Hence the first
term $a_0(t,q) c_{0,n-1}$ is in $\cI_n$.

Similarly, successively consider $\tilde \varphi_b$ with
$b=1,2,\dots, n-1$ we see that each term $a_i(t,q) c_{i,n-1}$ is in
$\cI_n$. \end{proof}

Let $\frac{\cI_n}{(q)_n}$ be the set of all rational function of the
form $f/(q)_n$, with $f\in \cI_n$. Let $\psi_d: \BZ[t^{\pm 1},
q^{\pm1}] \to \BZ[t^{\pm 1}]$ be the $\BZ[t^{\pm 1}]$-algebra
homomorphism defined by $\psi_d(q)=t^d$.

\begin{lemma} For every non-negative integer $n$ one has  $\psi_d(\frac{\cI_n}{(q)_n})\subset \Gd$.
\label{55}
\end{lemma}

\begin{proof} Fix a number $r$ co-prime with $d$. We need to show
that the multiplicity of $\Phi_r$ in the prime decomposition of the
denominator of
$$\psi_d \left(   \frac{c_{k,n}}{(q)_n}\right) = \frac{\prod_{i=0}^{k-1} (1-t^{id +1})}{\prod_{i=1}^{k} (1-t^{id})}$$
  is less than or equal to the that in the
numerator. Since $ 1- t^n = \prod_{m|n} \Phi_m$, we see that the
multiplicity of $\Phi_r$ in the denominator is the number of
elements of $\{d, 2d, \dots, kd \}$ which are divisible by $r$.
Since $r$ and $d$ are co-prime, this is the number of elements of
$\{1, 2, \dots, k \}$ which are divisible by $r$, and is equal to
$\lfloor k/r \rfloor$.

The multiplicity of $\Phi_r$ in the numerator is the number of
elements of $\{1, 1+d, \dots, 1+(k-1)d \}$ which are divisible by
$r$. This number is greater than or equal to $\lfloor k/r \rfloor$,
since any $r$ consecutive elements of any progressive sequence
contain one divisible by $r$.
\end{proof}

\begin{proof}(Proof of Proposition \ref{5})
By Theorem \ref{tech}, for every integer $b$,

$$ \tilde \varphi_b \frac{Z(k) (q)_k}{(q)_{2k+1}} \in \BZ[q^{\pm
1}].$$

It follows from  Proposition \ref{333} that $Z(k) (q)_k$ belongs to
$\cI_{2k+1}$, and hence $Y(k,d) \frac{(q)_k }{(q)_{2k+1}} = \psi_d
\left(\frac{Z(k) (q)_k}{(q)_{2k+1}}\right)$ belongs to $\Gd$, by
Lemma \ref{55}. It remains to notice that for every $n$, the
quotient $(q)_n/\{n\}!$ is a power of $q^{1/2}$, and
$\{2k+1\}!/\{k\}!$ is in $ q^{(k+1)(k+2)/4}\BZ[q]$.
\end{proof}

\newcommand{\Ld}{\Lambda_d}
\newcommand{\hLd}{\hat \Lambda_d}

\section{Cyclotomic completion}\label{500}
 After proving some additional facts in section
\ref{pre} we will  apply Habiro's result to our case.

\subsection{Cyclotomic completion, general results}
We fix a positive integer $d$. Recall that $\BN_d$ is the set of all
positive integers which are co-prime with $d$. Thus $\BN:=\BN_1$ is
the set of all positive integers. Recall also $R_d =\BZ[1/d]$.

In this section we will identify $q^{1/d}$ with $t$. Thus $A_d =
R_d[t^{\pm 1}]$, and $\Gd$ is obtained from $A_d$ by inverting all
the $\Phi_n(t), n\not \in \BN_d$. We also use $B_d:= R_d[t]$. One
has $B_d \subset A_d \subset \Gd$.

 For a subset $S\subset \BN$ let
$\Phi_S^*$ be the multiplicative set in $B_d$ generated by
$\Phi_r(t), r\in S$. Then $\Phi_S^*$ is a directed set with respect
to the divisibility relation $f|g$. For any ring $A$ containing
$\BZ[t]$, the principal ideals $(f)$ define a linear topology on
$A$.  Let

$$ A^S := \lim_{\overset{\longleftarrow }{f\in \Phi_S^*}} A/(f).$$

If $S'\subset S$ then there is a natural algebra homomorphism $
\rho_{S,S'}: A^S \to A^{S'}$. We want to know when $\rho_{S,S'}$ is
injective, for $A=B_d, A_d$, or $\Gd$.

\begin{proposition} (a) The map $\rho_{S, S\cap \BN_d}: \Gd^S \to \Gd^{S \cap
\BN_d}$ is an isomorphism.

(b) With $t=q^{1/d}$,
 $  \Gd^{\BN}=\Gd^{\BN_d}$ is equal to $\hGd$.

 (c) If $S\subset \BN_d$, then the inclusions $B_d \hookrightarrow
A_d \hookrightarrow \Gd$ induce  isomorphisms $B_d^S \cong A_d^S
\cong\Gd^S$.
 \label{1201}
\end{proposition}

\begin{remark} In particular $B_d^{\BN_d}= \hGd$. Note that in the
definition of $B_d^{\BN_d}$, we don't have any denominator. The ring
$B_d^S$ has been studied by Habiro \cite{Habiro2}. The author would
like to thank Habiro for pointing out the isomorphism
$B_d^{\BN_d}\cong \hGd$, which simplifies the original version of
the paper.
\end{remark}

\subsubsection{Injectivity of $\rho_{S,S'}$}

We say that $n,n'\in \BN$ are {\em p-adjacent} if $n/n'=p^e$, where
$p$ is a prime and $e$ an integer, and they are {\em adjacent} if
they are $p$-adjacent for some prime $p$.

For $S'\subset S \subset \BN_d$ we write $S'\prec S$ if for every
$n\in S$ there are $n_1,\dots,n_k \in S$ and $n_{k+1}\in S'$, such
that $n_i$ and $n_{i+1}$ are adjacent, and $n_1=n$.

For $m,n \in \BN$ and any ring $R \supset \BZ$, Habiro
\cite{Habiro2} defined $m \Leftrightarrow_{R} n$ if $m,n$ are
$p$-adjacent and $R$ is $p$-adically separated, i.e. $\cap_{j\ge 0}
p^j R$ consists of 0 only. It is easy to see that if $p$ is a prime,
$p \in \BN_d$, then $R_d=\BZ[1/d]$ is $p$-separated. Hence if $m, n$
are in $\BN_d$ and they are adjacent, then $m \Leftrightarrow_{R_d}
n$ in Habiro's sense.

Suppose that $S_0\prec S\subset \BN_d$. From theorem 4.2 of
\cite{Habiro2} one has $\rho_{S,S_0}: B_d^S \to B_d^{S_0}$ is
injective. Using Proposition \ref{1201} part (c) we get the
following.

\begin{theorem} Suppose that $S_0\prec S\subset \BN_d$.  Then
  $\rho_{S,S_0}: \Gd^S \to \Gd^{S_0}$
is injective. \label{73}
\end{theorem}

\subsubsection{Evaluation map}

\newcommand{\tev}{\tilde \ev}

Suppose  $\xi$ is a root of order $r$, which belongs to $S \subset
\BN_d$. By Proposition \ref{1} below, $\Gd/(\Phi_r) \cong
\BZ[1/d][t^{\pm 1}]/(\Phi_r) $. The last one is isomorphic to
$\BZ[1/d][\xi]$ via $t \to \xi$. We can define the evaluation map
$\tev_\xi: \Gd^S \to \BZ[1/d][\xi]$ by the composition

$$ \Gd^S \to \Gd^{(\Phi_r)} \to \Gd/(\Phi_r) \equiv
\BZ[1/d][\xi].$$

Since $r$ and $d$ are co-prime, $\zeta:= \xi^d$ is also a root of
order $r$, and by definition, $\tev_\xi(f) = \ev_\zeta(g)$, where
$g$ is obtained from $f$ by the substitution $t\to q^{1/d}$. Hence
we get the following

\begin{proposition} For $g\in \hGd$ and $\xi\in U_d$, one has
$\ev_\xi(g) \in \BZ[1/d][\xi]$.
\end{proposition}

Suppose $\Omega$ is a set of roots of unity whose orders form a set
$T$, which is a subset of  $S \subset \BN_d$. Using the evaluation
at every element in $\Omega$ we can define
$$\tev_{\Omega}:\Gd^S \to
\prod_{\xi\in \Omega}\BZ[1/d][\xi], \qquad \tev_{\Omega}(f) =
(\tev_\xi(f), \xi \in \Omega).$$

Again using $B_d^S \cong \Gd^S$, Theorem 6.1 of \cite{Habiro2} gives
us the following.

\begin{theorem} Suppose $\Omega,T,S$ are as above. Assume that there
is $n\in S$ such that $\{n\} \prec S$ and there are infinitely many
elements in $T$ adjacent to $n$.  Then $\tev_{\Omega}: \Gd^S \to
\prod_{\xi\in \Omega}\BZ[1/d][\xi]$ is injective.\label{last}
\end{theorem}

Theorem \ref{main3} is a special case of this theorem with $S=
\BN_d, n=1$.

\subsection{Proof of Proposition \ref{1201}} \label{pre}
Part (a). Since $\Phi_r(t)$ is invertible in $\Gd$ if $r$ is not in
$\BN_d$, we have $ \Gd^S = \Gd^{S \cap \BN_d}$.

Part (b).  Since $ 1-t^n = \prod_{r|n}\Phi_r(t)$, one has
$\frac{1-t^{nd}}{1-t^{n}}$ is the product of several $\Phi_r$ with $
r \not \in \BN_d$. Hence $\frac{1-t^{nd}}{1-t^{n}}$ is invertible in
$\Gd$. This means, if $q=t^d$ then the elements $(t;t)_n$ and
$(q;q)_n$ defines the same principal ideal in $\Gd$, and hence they
define the same completion, or $\Gd^{\BN}=\hGd$.

 The rest of this subsection is devoted to a proof of part
(c).

\begin{lemma}\label{90}

(a) If $m,n$ are not adjacent, then $(\Phi_m) + (\Phi_n) =(1)$ in
$\BZ[t]$.

(b) If $m=np^e$, where $e>0$ and $p$ a prime, then $(\Phi_m) +
(\Phi_n) =(p) + (\Phi_n)$ in $\BZ[t]$.
\end{lemma}

\begin{proof} (a) is a well-known fact.

(b) Let $a= t^{n p^{e-1}}$, then $\Phi_m$ is a divisor of
$$ g(t) = \frac{1-a^p}{1-a}= 1+a+\dots+ a^{p-1}.$$

Since $a(\xi)=1$, and hence $g(\xi)=p$ if $\xi$ is a root of unity
of order $n$, we have that $g(t) =p \pmod {\Phi_n}$.  This means
$\Phi_m$ divides $p$ in $\BZ[t]/(\Phi_n)$. It is known that there is
a positive $k$ such that  $\Phi_m = \Phi_n^k  \pmod{p}$. It follows
that $p$ divides $\Phi_m$ in $\BZ[t]/(\Phi_n)$. Thus in $\BZ[t]/(
\Phi_n)$, $\Phi_m$ and $p$ define the same principal ideal. This is
equivalent to $(\Phi_m) + (\Phi_n) =(p) + (\Phi_n)$ in $\BZ[t]$.
\end{proof}

\begin{corollary} Suppose $m_1,\dots,m_k$ are co-prime with $d$ and $n$ is not, i.e. $m_1,\dots,m_k\in \BN_d$ and $n\not \in \BN_d$.
Then $\Phi_n$ is invertible in $\BZ[1/d][t]/(f)$, where $f=
\prod_{i=1}^k \Phi_{m_i}$. Here $m_1,\dots,m_k$ are not necessarily
distinct. \label{1200}
\end{corollary}

\begin{proof} Note that in a commutative ring $A$, an element $a$ is invertible
in $A/(b)$ if and only if   $(a)+(b)=(1)$.  If $(a) + (b)=(1)$ and
$(a)+(c)=(1)$ then, multiplying together, one gets $(a) + (bc)
=(1)$. Hence it's enough to consider the case $k=1$, with $m_1=m$.

If $m,n$ are not adjacent, then  $(\Phi_n) + (\Phi_r) =(1)$ in
$\BZ[t]$, hence we are done in this case.

Suppose $n/m = p^e$, with $p$ a prime.  Then $e>0$, since otherwise
$n$ is a divisor of $m$ and hence  is co-prime with $d$. By the same
reason, $p$ is a divisor of $d$.  By Lemma \ref{90}, $(\Phi_m) +
(\Phi_n) \supset (p)$ in $\BZ[t]$. Since $p$ divides $d$, it is
invertible in $\BZ[1/d][t]$, hence $(\Phi_m) + (\Phi_n) =(1)$ in
$\BZ[1/d][t]$.
\end{proof}

\begin{proposition}\label{1}
Suppose $f= \prod_{i=1}^k \Phi_{m_i}$ with $m_i\in \BN_d$. Then the
homomorphism  $A_d/(f) \to \Gd/(f)$, induced from the embedding
$A_d\hookrightarrow \Gd$, is an isomorphism.
\end{proposition}

\begin{proof} Recall that $\Gd$ is obtained from $A_d$ by inverting all the $\Phi_n, n\not \in
\BN_d$.

By Corollary \ref{1200}, every $\Phi_n, n\not \in  \BN_d$, is
invertible in $A_d/(f)$. This proves the surjectivity. Injectivity
follows easily from the fact that every $\Phi_n, n\not \in  \BN_d$,
is co-prime with $f$.
\end{proof}

\begin{corollary} If $S\subset \BN_d$, then the inclusion $A_d \hookrightarrow
 \Gd$ induces an isomorphism $A_d^S \cong \Gd^S$.\label{20200}
\end{corollary}

\begin{proposition} For any  $S\subset \BN$, the inclusion $B_d \hookrightarrow
 A_d$ induces an isomorphism $B_d^S \cong A_d^S$.\label{20201}
\end{proposition}

\begin{proof} Note that $\Phi_n(0)=\pm 1$, one has $(t) +
(\Phi_n(t))=(1)$ in $\BZ[t]$. It follows that if $f$ is the product
of several $\Phi_n(t)$, then $t$ is invertible in $\BZ[t]/(f)$, and
hence in $B_d[t]/(f)$. It follows that $B_d^S \cong A_d^S$.
\end{proof}

Part (c) of Proposition \ref{1201} follows from Corollary
\ref{20200} and Proposition \ref{20201}.

\subsection{Proof of Theorem \ref{main4}} Part (a).  Using Proposition \ref{1201}, with $t=q^{1/d}$, we can
identify  $\hGd$ with $\Gd^{\BN_d}$. The Taylor expansion map $T_1$
is then the map $\rho_{\BN_d, \{1\}}: \Gd^{\BN_d} \to \Gd^{\{1\}}$.
The image space $\Gd^{\{1\}}$ is equal to $\BZ[1/d][t^{\pm
1}]^{(\Phi_1)} = \BZ[1/d][[t-1]]$, by Proposition \ref{1}.

As explained, to convert power series in $t-1= q^{1/d}-1$ into power
series in $q-1$, one use

$$ q^{1/d}-1 = (1+(q-1))^{1/d} -1 = \sum_{n=1}^\infty
\binom{1/d}{n} (q-1)^n.$$

It is well-known that $\binom{1/d}{n} \in \BZ[1/d]$ for every $n$.
Hence the image of $T_1$ is in $\BZ[1/d][[q-1]]$.

Part (b), the most difficult part, is a special case of Theorem
\ref{73}, with $S=\BN_d$ and $S'=\{1\}$.

Part (c) follows from the definition of the Ohtsuki series
\cite{Ohtsuki1}. \qed

\ifx\undefined\bysame
    \newcommand{\bysame}{\leavevmode\hbox
to3em{\hrulefill}\,}
\fi

\end{document}